\documentclass[11pt,a4paper]{article}
\usepackage[T1]{fontenc}
\usepackage[utf8]{inputenc}
\usepackage[english]{babel}
\usepackage{lmodern,microtype}
\usepackage{amsmath,amssymb,amsthm,mathtools}
\usepackage[margin=25mm]{geometry}
\usepackage{graphicx,booktabs,array,tabularx,longtable,enumitem,placeins}
\usepackage[numbers,sort&compress]{natbib}
\usepackage{xcolor}
\usepackage{hyperref}
\hypersetup{hidelinks,
 pdftitle={Balanced cut trees for classification with local neural networks},pdfauthor={Oleg R. Musin}}
\setlist{itemsep=2pt,topsep=4pt}
\newtheorem{theorem}{Theorem}[section]
\newtheorem{proposition}[theorem]{Proposition}

\theoremstyle{definition}
\newtheorem{definition}[theorem]{Definition}

\theoremstyle{remark}

\newcommand{\R}{\mathbb R}

\DeclareMathOperator*{\argmax}{arg\,max}

\DeclareMathOperator{\dist}{dist}

\newcommand{\relu}{\operatorname{ReLU}}
\title{Balanced cut trees for classification with local neural networks}
\author{Oleg R. Musin}
\date{}

\begin{document}
\maketitle
\begin{abstract}
We study binary trees in which each cut bisects every class of a finite
training sample. Each leaf therefore retains the original class proportions.
A new observation is classified by a single local network, using its
coordinates within the leaf. We give conditions for strict polynomial
bisection and show that a one-hidden-layer ReLU function of width at most
$k$ strictly bisects $k$ even classes of points with distinct coordinates.
Gaussian smoothing leads to a continuation method for proposing cuts;
a finite test determines whether a prescribed coefficient box contains
an exact bisector. For classes lying near separated smooth graphs, we
give a bound on the size of a region in which affine scores suffice.
Experiments illustrate both the benefit and the cost of requiring
balance. Small local models can classify geometrically distinct regions
at low prediction cost, but exact bisection may split those same regions
when their empirical class counts vary.
\end{abstract}

\section{Introduction}\label{sec:introduction}

\subsection{The ham-sandwich theorem and finite data}

Suppose that observations in $\R^d$ belong to $k$ classes, represented by
finite sets $X_1,\ldots,X_k$. A classifier assigns a class to a new point
from its coordinates. The usual approach is to fit one function to all
the observations. We consider a geometric alternative: divide the data
into regions, fit a simpler classifier in each region, and use the
classifier belonging to the region of the new point. The question is
how to choose the regions.

The ham-sandwich theorem suggests a particularly symmetric choice.
Instead of separating the classes at an internal node, we divide
\emph{each} class into two equal parts. Repeating the construction gives
regions with the same number of observations of each fixed class at
a given depth. Classification takes place within these regions. The
purpose of the division is therefore to simplify the local geometry
while retaining observations from every class.

The starting point is the following classical theorem; see
\cite{StoneTukey1942,Steinhaus1945,Matousek2003}.

\begin{theorem}[Ham-sandwich theorem]\label{thm:hs-intro}
Let $\mu_1,\ldots,\mu_k$ be finite absolutely continuous Borel measures
on $\R^d$, where $k\le d$. There is an affine hyperplane whose two open
halfspaces have equal $\mu_i$-measure for every $i$.
\end{theorem}

The short topological proof also indicates why nonlinear versions are
possible. Represent an oriented affine hyperplane by
$a=(a_0,a_1,\ldots,a_d)\in S^d$, and consider the difference between the
two masses determined by the sign of
$a_0+a_1x_1+\cdots+a_dx_d$. These $k$ differences form a continuous odd
map from $S^d$ to $\R^k$. Continuity follows because a hyperplane has
measure zero. The Borsuk--Ulam theorem gives a zero of this map, after
adding zero coordinates when $k<d$. If at least one measure is nonzero,
the two constant functions at the poles are not zeros, so the resulting
function defines a genuine hyperplane. If all measures vanish, any
hyperplane suffices.

For point sets, the corresponding statement allows points on the
hyperplane: each open halfspace contains at most
$\lfloor |X_i|/2\rfloor$ points of $X_i$. This familiar discrete version
does not by itself assign the boundary points to two equal groups. For
classification we need an actual partition of the original
observations, determined by their coordinates. Even cardinalities and
appropriate general position allow a strict division with no points on
the boundary. Repeated coordinates can prevent it. We state the finite
conditions explicitly in Section~\ref{sec:models}.

This distinction is useful even in the simplest examples. A line may
bisect two planar classes, but three classes generally require a more
flexible curve. On the other hand, a complicated curve chosen solely to
give equal counts may divide a geometrically simple class into many
pieces. Existence of a balanced division is the first part of the
problem. Choosing a division that helps classification is the second.

\subsection{Balanced trees and local classification}

The classifier consists of a binary tree of cuts and a collection of local
networks, one for each leaf. At a branch, the sign of a function of the
coordinates determines which of two regions contains the point. We choose
this function so that each region receives half the training observations
of each class present at the branch. At a leaf, a network fitted to the
local data takes the original coordinates as input and gives one score
for each class. The class with the highest score is the prediction.
Thus evaluating the classifier at a new point requires only the cuts
along its path and the network at its leaf.

The depth can depend on the region. If a local network already classifies
a region with sufficiently small error, we may keep it as a leaf and
continue dividing elsewhere. A leaf reached after $q$ strict bisections
contains $|X_i|/2^q$ observations of class $i$.

Why might this help? Consider two separated curves whose relative
position changes across the plane. A single classifier must describe
the whole boundary between them. In a small neighbourhood, the two
curves may be well approximated by straight lines, and their separation
can be decided by a much simpler function. The same idea applies to
several nearby surfaces. The advantage comes from allowing different
simple rules in different regions and evaluating only the rule that is
needed at the new point.

Balance gives each local problem a prescribed share of every class.
It does not imply that a region is small or that its class boundaries
are simple. Indeed, the class proportions in each leaf are exactly the
same as in the full sample. Consequently the identity of the leaf alone
provides no information about the class on that sample. This is a
fundamental difference from a decision tree whose leaves are intended
to contain mostly one class. Here the useful information remains in
the coordinates within the leaf.

For two classes we suggest the following working hypothesis. Balanced
cuts are especially promising when the classes lie close together over
a substantial fraction of the mass of each, the direction of separation
changes repeatedly across the data, and admissible cuts isolate regions
with simple local decision rules. Establishing when such a partition
offers an advantage over a global classifier will require both extensive
experiments and a mathematical theory.

Two choices must therefore be made together: the family of allowed
cuts and the family of local classifiers. Regions appropriate for
affine separation need not be best for small nonlinear networks. A
tree with many leaves also stores many local models, even though one
prediction uses only one of them. Thus prediction cost, storage, and
the effort required to construct the tree are separate quantities.

Dividing a statistical problem among local models has a long history.
Mixtures of local experts and their hierarchical versions were studied
in~\cite{JacobsEtAl1991,JordanJacobs1994}; affine tests in decision trees
are discussed in~\cite{MurthyKasifSalzberg1994}. Soft decision trees use
weighted combinations of predictions from different branches
\cite{IrsoyYildizAlpaydin2012,FrosstHinton2017}. Snopov and Musin~\cite{SnopovMusin2025Activations} introduced activation
functions that split data representations within a neural network.
Our earlier work
\cite{MusinSnopov2025Decomposition} considered geometric decomposition
for neural classification. The present paper studies the additional
requirement of simultaneous class bisection and its consequences for
the construction, local training, and cost of prediction.

There is a further practical reason to distinguish a partition from
the data used to fit its local classifiers. A classifier fitted only
to points on one side of a cut may be poorly determined near the
boundary. We therefore allow it to use nearby observations from the
adjacent region. The regions themselves remain disjoint, and the new
point still uses a single local classifier. Such overlapping training
sets are a natural way to reduce a boundary effect. They are helpful
when neighbouring observations follow a compatible local rule; they
can be harmful when they impose a different rule. Section~\ref{sec:models}
makes this distinction precise.

\subsection{Nonlinear bisectors}

Polynomial cuts extend the geometric construction to more classes.
For example, the map
\[
                   (x,y)\longmapsto (x,y,x^2+y^2)
\]
turns a circle or a line in the plane into the inverse image of a plane
in three dimensions. Simultaneous circular bisection of three planar
sets, with a straight line allowed as a limiting case, already appears
in Steinhaus's work~\cite{Steinhaus1945} in 1945.

More generally, a polynomial of degree at most $s$ in $d$ variables
has $\binom{d+s}{s}$ coefficients. The polynomial ham-sandwich theorem
provides a simultaneous bisector for $k$ absolutely continuous measures
when
\[
                        \binom{d+s}{s}-1\ge k.
\]
Its proof uses the same odd-map argument as Theorem~\ref{thm:hs-intro},
now on the sphere of polynomial coefficients. A nonzero polynomial
has a zero set of Lebesgue measure zero. This observation establishes
continuity in the original space; the lifted measures need not be
absolutely continuous in the larger space.
Polynomial and function-space partitions are classical tools in
geometry~\cite{StoneTukey1942,Golasinski2006,KarasevHubardAronov2014,
GuthKatz2015,KaplanMatousekSharir2012}. Their role here is to provide
families of balanced regions for local classification.

Neural cuts have a similar interpretation. A hidden layer gives a
collection of functions of the original coordinates, and an affine
combination of these functions defines a cut. For the rectified linear
function $\relu(t)=\max\{t,0\}$, which takes the positive part of a real input $t$, the classical theorem of Hobby and
Rice~\cite{HobbyRice1965} leads to a particularly simple conclusion:
a suitably chosen one-hidden-layer cut of width at most $k$ can bisect
$k$ absolutely continuous measures. There is also a strict finite
version for even classes whose combined observations are pairwise
distinct. We give the construction, including its finite interpretation,
in Section~\ref{sec:models}.

The dependence on the data is essential. A prescribed collection of
features may identify distinct observations and make exact bisection
impossible. Conversely, the existence of an adaptive cut says little
about the geometry of the resulting regions. These are reasons to
study the choice of features and cuts as part of the classification
problem. In the main comparisons below, both the global network and
the local networks receive the original coordinates; any additional
functions used to determine the cuts are included in the cost of the
tree.

\subsection{Construction and scope of the paper}

Topological existence theorems leave a computational problem.
For two planar point sets, the exact algorithm of Lo, Matou\v{s}ek,
and Steiger~\cite{LMS1994} has linear complexity. For $d>2$, their running-time bound depends on the cost of
constructing a level in an arrangement of hyperplanes in dimension
$d-1$. The difficulty already arises for $k>3$ classes. To apply the
ham-sandwich theorem to $k$ classes in general, the feature space
must have dimension at least $k$. Thus even for four planar classes
the lifted hyperplane search takes place in at least four dimensions.
Relying on the exact algorithm alone can already make the search for
a common cut difficult in this case.

The author's \emph{Finding solutions with the Borsuk--Ulam theorem}~\cite{MusinFindingSolutions}
provides a numerical alternative, building on the geometric approach
of~\cite{MusinManifolds2012}. Its starting point
is a reference problem with a known regular solution. A deformation
joins the reference balance equations to those for the desired cut;
following their common zeros gives a numerical candidate. For point
data the counts jump whenever a cut crosses an observation. Replacing
each lifted point by a Gaussian makes the equations smooth enough to
follow, and the signs at the original points decide whether the
candidate actually bisects each class. Section~\ref{sec:certified-construction}
gives the equations and a finite search when a candidate needs adjustment.

A scan of the data evaluates the balance equations at a given cut, but
the number of continuation steps depends on the path, and a poorly
conditioned correction can require a smaller step. For the general
discrete ham-sandwich problem with dimension part of the input,
Filos-Ratsikas and Goldberg proved PPA-completeness~\cite{FRG}.
This worst-case result concerns the discrete search problem; the work
of following a particular smoothed path is a separate question.

Gaussian smoothing can also describe the intended statistical problem.
For instance, a measurement may represent an uncertain location rather
than an exact point. Then dividing a specified distribution of possible
locations may be more meaningful than dividing the observed centres.
In that setting a fixed smoothed measure is the object to be bisected.
An exact division of its centres need not improve its mass balance.
We discuss both uses of smoothing while keeping their objectives clear.

Section~\ref{sec:models} develops the mathematical construction. We
begin with the finite partition and the elementary consequence for class
proportions, then describe polynomial and neural bisectors. We relate
smooth balance equations to finite counts and explain how a bounded
search can verify or exclude a strict cut in a prescribed coefficient
region. A local separation result for smooth graphs explains one way
in which smaller regions simplify classification. We then discuss
overlap, the learning of cuts, and computational cost. The classical
existence theorems are used as foundations; the point of the analysis
is to connect them to the requirements of a local classifier.

Section~\ref{sec:experiments}, \emph{Experiments and future work},
describes a selection of examples. They include planar curves,
several classes with different local decision rules, algebraic cuts,
and observed data. Some examples show that small local networks can
give an economical prediction rule. Others show how the apparent
advantage changes with a stronger global network, a new sample, or a
different choice of cuts. These comparisons motivate the central
remaining questions: how to find useful balanced partitions reliably,
how to learn them together with local classifiers, and how to evaluate
their advantages on independent observations.

The experiments are exploratory. Their purpose is to identify useful
geometric mechanisms and concrete difficulties. A general superiority
claim would require substantially broader evidence. We retain enough
detail to specify the models and interpret the comparisons, and leave
the publication of a complete collection of programs and data to a
separate stage of the project.

\section{Balanced partitions and local classifiers}\label{sec:models}

\subsection{The finite sample and its partition}

Let $X_1,\ldots,X_k\subset\R^d$ be the nonempty classes in a finite training
sample, and put $n_i=|X_i|$ and $N=\sum_i n_i$.
We regard the observations as indexed points, so that repeated coordinates
are counted with their multiplicities. The label of a point records the
class to which it belongs. Labels are known during construction; the
classifier must predict the label of a new point from its coordinates.

A rooted binary tree divides the input space recursively. At an internal
node $v$ choose a real-valued function $g_v$ and send a point to the left
when $g_v(x)\le0$, and to the right when $g_v(x)>0$. The functions may be
affine, polynomial, or represented by a small neural network. Their signs
determine disjoint regions $C_v$ at the leaves, whose union is $\R^d$.
Write $v(x)$ for the leaf containing $x$.

At each leaf fit a function $h_v:\R^d\to\R^k$. Its coordinates are scores
for the $k$ classes: the predicted class is the coordinate with the largest
score, with a fixed order used to resolve ties. The resulting classifier is
\begin{equation}\label{eq:cut-network}
 F(x)=h_{v(x)}(x),\qquad
 \widehat y(x)=\argmax_{1\le i\le k}F_i(x).
\end{equation}
Thus a query uses one path through the tree and one local classifier.
An affine choice $h_v(x)=W_v^{\mathsf T}x+b_v$ gives the simplest local
model. For a nonlinear local model, consider a neuron with weight vector
$w\in\R^d$ and bias $b\in\R$. It first forms the real number
$t=w^{\mathsf T}x+b$ and then outputs its positive part,
$\relu(t)=\max\{t,0\}$. This operation is called a rectified linear
unit (ReLU). A hidden layer of width $m$ consists of $m$ such neurons,
each with its own $w$ and $b$. For example, a local network with one
hidden layer gives the vector of class scores
\[
 h_v(x)=c_v+\sum_{j=1}^{m}a_{vj}\,
           \relu(w_{vj}^{\mathsf T}x+b_{vj}),\qquad a_{vj},c_v\in\R^k.
\]
Each leaf has its own fitted coefficients.
Throughout the principal construction these functions receive the
original coordinates $x$, just as the global network used for comparison
does. Features introduced to construct a cut need not also be inputs to
the local classifiers.

\begin{definition}\label{def:balanced-tree}
A finite tree is \emph{strictly balanced} if, at every internal node,
the cut is nonzero at all training points reaching that node and divides
each of their classes into equal halves. Equivalently, if $X_{v,i}$ is
the class-$i$ sample reaching $v$, then
\begin{equation}\label{eq:strict-balance}
 g_v(x)\ne0\quad(x\in X_{v,i}),\qquad
 |X_{v0,i}|=|X_{v1,i}|=\frac{|X_{v,i}|}{2}.
\end{equation}
The children are formed by the signs of the same function $g_v$ for every
class.
\end{definition}

The last requirement matters. Assigning points to two lists with the
desired class counts is easy, but such lists need not be separated by an
admissible function. The construction must supply a geometric rule that
also applies to new points. Identical coordinates therefore follow the
same path, and the label is never used to choose a branch.

Induction gives a useful elementary consequence of the definition.
A node at depth $q$ contains exactly $n_i/2^q$ training points of class
$i$. A complete tree of depth $H$ consequently requires $2^H\mid n_i$
for every $i$. The classes themselves need not have the same size.
Every node preserves their original proportions $n_i/N$.
Indeed, for a leaf $v$ at depth $q(v)$,
\[
 \frac{|X_{v,i}|}{|X_i|}=2^{-q(v)},\qquad
 \frac{|X_{v,i}|}{\sum_j|X_{v,j}|}=\frac{n_i}{N}.
\]
These identities apply separately to each leaf, even when the leaves
have different depths. For example, during construction one may fit a
local network at a candidate node and estimate its error by the fraction
of misclassified observations among validation points reaching that node.
If this error is below a prescribed tolerance, the node may be declared
a leaf. Other branches may still be divided. With a depth limit $H$,
the resulting tree therefore has at most $2^H$ leaves. Strict balance
is required only at the cuts actually made.

This explains the role of balance in the present method. In an ordinary
classification tree, one seeks regions dominated by a single class.
Here every region retains all the original class proportions, and knowing
the leaf alone gives no class information on the training sample. The
possible gain comes from using the coordinates \emph{within} a region:
a complicated global decision boundary may consist of simpler local
pieces. Balance supplies data from each class to each local problem,
while the choice of the regions determines whether those problems become
easier. The two issues must be considered together.

\begin{figure}[htbp]
\centering
\includegraphics[width=\linewidth]{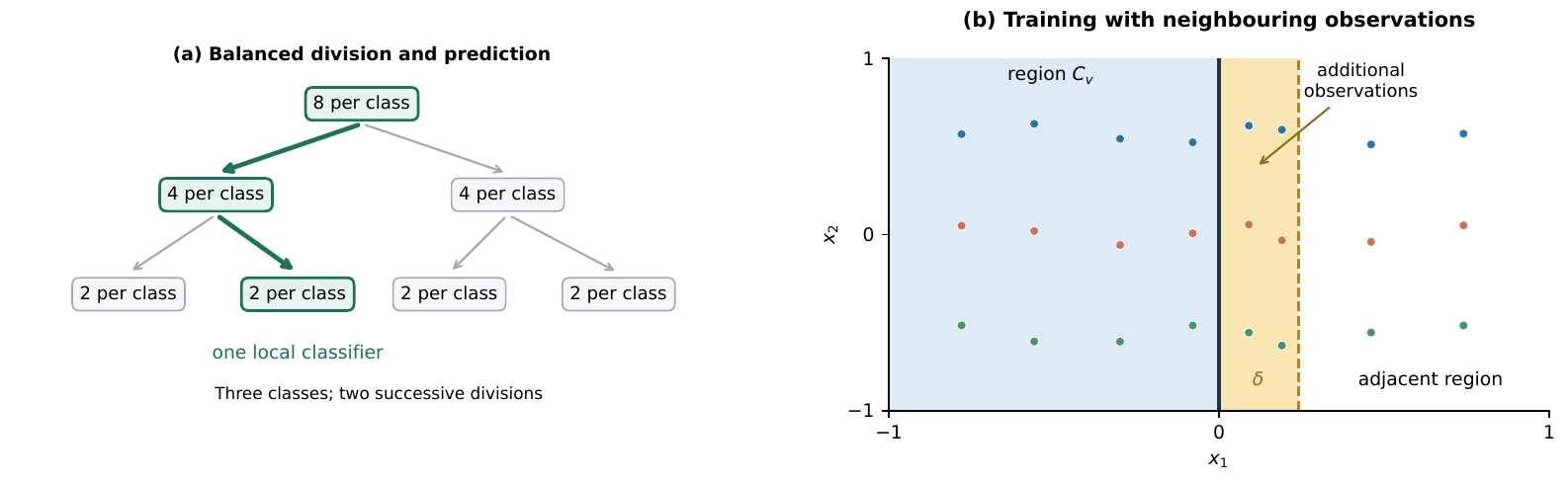}
\caption{A balanced tree and local fitting with overlap. A query follows
one path and is classified in its final region. Nearby training points
from an adjacent region may help fit the local classifier, while the
partition used for prediction remains fixed.}
\label{fig:routing-overlap}
\end{figure}

\subsection{Polynomial cuts and exact finite bisection}\label{sec:gate-families}

Fix a family of real functions
$\phi=(\phi_1,\ldots,\phi_M):\R^d\to\R^M$ and consider cuts
\[
 g(x)=b+a^{\mathsf T}\phi(x).
\]
The map $\phi$ changes the functions available for making a partition.
It need not be an embedding. With $\phi(x)=x$ the cuts are hyperplanes.
In the plane, the map $(x,y)\mapsto(x,y,x^2+y^2)$ also permits circles.
This is the lifting interpretation of the classical circular bisection
theorem discussed in the introduction.

More generally, one fixed nonzero quadratic form $Q$ gives
\begin{equation}\label{eq:shared-lift}
 \phi_Q(x)=(x,x^{\mathsf T}Qx),\qquad
 g_v(x)=b_v+w_v^{\mathsf T}x+\alpha_vx^{\mathsf T}Qx.
\end{equation}
Using the same $Q$ at every node restricts the quadratic parts to be
proportional. Allowing a separate form $Q_v$ at each node is a larger
family. For example, independent planar quadratic cuts use the five
nonconstant features $(x,y,x^2,xy,y^2)$.

The natural general setting is a polynomial space
\[
 \mathcal P=\operatorname{span}\{1,\phi_1,\ldots,\phi_M\}
\]
of dimension $M+1$. The number of nonconstant coefficients determines
how many measures can be bisected simultaneously. Some familiar choices
are recorded in Table~\ref{tab:dictionaries}.

\begin{table}[htbp]
\centering
\begin{tabular}{lll}
\toprule
Functions used for cuts & $M$ & Number of measures\\
\midrule
Affine functions & $d$ & $k\le d$\\
Affine functions and one quadratic form & $d+1$ & $k\le d+1$\\
All quadratics & $\binom{d+2}{2}-1$ & $k\le M$\\
All polynomials of degree at most $s$ & $\binom{d+s}{s}-1$ & $k\le M$\\
\bottomrule
\end{tabular}
\caption{Polynomial spaces for simultaneous bisection. The final column
concerns absolutely continuous measures. Finite samples require the
additional conditions in Theorem~\ref{thm:finite-polynomial}.}
\label{tab:dictionaries}
\end{table}

\begin{theorem}[Polynomial balanced partitions]\label{prop:arbitrary-family-lift}
Let $\mu_1,\ldots,\mu_k$ be finite positive absolutely continuous Borel
measures on $\R^d$, and let $\mathcal P$ be a polynomial space of
dimension $M+1$ containing the constants, with $k\le M$.
For every integer $H\ge0$ there is a complete tree of depth $H$ whose
cuts belong to $\mathcal P$ and whose leaves satisfy
\[
 \mu_i(C_v)=2^{-H}\mu_i(\R^d),\qquad 1\le i\le k.
\]
\end{theorem}

\begin{proof}
A nonzero polynomial has a Lebesgue-null zero set. Hence the differences
between the masses of its positive and negative regions depend
continuously on its coefficient vector on the unit sphere $S^M$.
Changing the sign of the coefficients reverses every difference.
The Borsuk--Ulam argument therefore supplies a common bisector; this is
the polynomial ham-sandwich theorem
\cite{StoneTukey1942,Matousek2003}.
A nonzero constant has the same strict sign everywhere and cannot bisect a measure of positive mass.
The restriction of each $\mu_i$ to either resulting region is again
absolutely continuous. Repeating the same argument at every node proves
the assertion by induction. Polynomial boundaries have zero measure,
so assigning them to the left child does not change the masses.
\end{proof}

The argument takes place in the original space $\R^d$. The images of
the measures under $\phi$ are often singular in $\R^M$; absolute
continuity in the larger space is not needed. No convexity or
connectedness assumption is imposed on the original supports or the
regions of the tree.

For a finite sample, points lying on the cut cannot simply be ignored.
The discrete ham-sandwich theorem allows such points and bounds the
number on either open side by one half of the class size. The following
condition makes it possible to move the cut slightly and assign every
sample point to a strict side.

\begin{theorem}[Strict finite bisection in a feature space]\label{thm:finite-polynomial}
Let $\mathcal P=\operatorname{span}\{1,\phi_1,\ldots,\phi_M\}$ have
dimension $M+1$, with $k\le M$. Suppose that every $n_i$ is positive
and even, and that every collection of at most $M+1$ evaluation rows
\[
 e(x)=(1,\phi_1(x),\ldots,\phi_M(x))
\]
from the whole indexed sample is linearly independent. The features may be arbitrary real-valued functions. Then a function
in $\mathcal P$ strictly bisects every class. If $2^H\mid n_i$ for
all $i$, the construction can be repeated to give a complete strictly
balanced tree of depth $H$.
\end{theorem}

\begin{proof}
Pairwise independence makes the lifted points distinct. If $k<M$,
adjoin $M-k$ singleton classes at distinct auxiliary points and discard
their constraints after applying the theorem.
Apply the discrete ham-sandwich theorem to the lifted points
\cite[Section~3.1]{Matousek2003}. It gives a nonzero $p\in\mathcal P$
with at most $n_i/2$ points of class $i$ on either open side.
The independence assumption implies that at most $M$ sample points
lie on $p=0$. If class $i$ has $a_i$ positive and $b_i$ negative
values, assign positive signs to $n_i/2-a_i$ of its zero values and
negative signs to the rest. Both numbers are admissible.

The evaluation rows of the boundary points are independent. There is
therefore $u\in\mathcal P$ taking the prescribed values $+1$ and $-1$
at those points. For sufficiently small $\varepsilon>0$, the function
$p+\varepsilon u$ has those signs on the boundary and preserves all
previously nonzero signs. It is the required strict bisector.
Every subset of the sample inherits the independence assumption, so the
same proof applies recursively.
\end{proof}

For a \emph{fixed} polynomial space, the independence assumption holds
with probability one if the joint distribution of the sample coordinates
has a density. To see this, note that the vectors $e(x)$ span
$\R^{M+1}$; otherwise a nonzero polynomial in $\mathcal P$ would vanish
identically. For any prescribed number of rows at most $M+1$, some
minor is thus a nonzero polynomial in their coordinates. Its zero set
has measure zero, and there are only finitely many subsets of the
sample to consider. A feature family chosen after seeing the sample
requires its own independence argument or a direct check of the
resulting signs.

\begin{figure}[htbp]
\centering
\includegraphics[width=\linewidth]{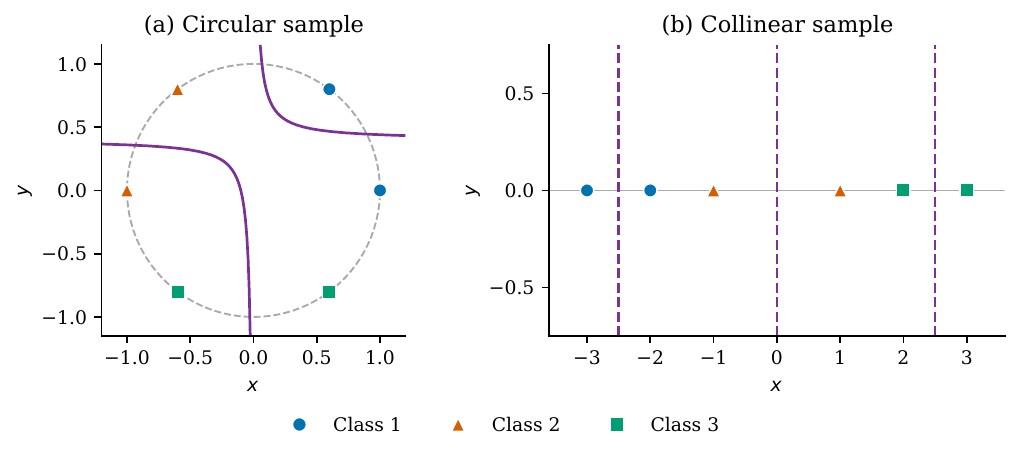}
\caption{The feature family matters for finite samples. Left: the three
pairs on a circle have a common hyperbolic bisector, but no common
circular or linear bisector. Right: three separated pairs on a line
require a polynomial of degree at least three.}
\label{fig:finite-lift-obstructions}
\end{figure}

The finite assumptions have a geometric meaning. On the unit circle,
the feature $x^2+y^2$ is constant; circular cuts then act on the sample
like affine cuts. The six points in Figure~\ref{fig:finite-lift-obstructions},
labelled $A,A,B,B,C,C$ in cyclic order, have no strict circular or linear
bisector of all three pairs. Any half-plane selects a consecutive block
of the cyclic order, whereas no block of length three contains one point
from every pair. A hyperbola, using the feature $xy$, does bisect them.
Likewise, three separated pairs on a line require at least three roots
along that line, so a quadratic cannot strictly bisect them all.
These examples explain why the choice of polynomial space matters for
finite data supported on lower-dimensional sets.

Repeated coordinates give a different obstruction, independent of
polynomial degree. All copies of one coordinate must remain together.
For example, suppose three distinct coordinates have class multiplicity
vectors $(2,1)$, $(1,2)$ and $(1,1)$. Both total class counts are four,
but no subset of these three whole groups has counts $(2,2)$. No
deterministic cut can bisect the sample. Thus divisibility is necessary
but is not sufficient when identical inputs are present.

\subsection{Small neural networks as common bisectors}\label{sec:relu-lifting}

The same neurons can be used at an internal node. Here their outputs
are combined into a single real number whose sign selects a branch.
In the theorem below, $t=u^{\mathsf T}x$ denotes the coordinate obtained
by projecting the observation onto a fixed direction. A neuron of the
form $\relu(t-s)$ becomes active when this coordinate exceeds the
threshold $s$. An affine combination of features
$\relu(w_j^{\mathsf T}x+b_j)$ is a one-hidden-layer ReLU function.
Consequently a hidden layer is another way of choosing a family of
features for the cuts. The next result connects this elementary
observation with one-dimensional simultaneous bisection.

\begin{theorem}[A ReLU bisector of width at most $k$]
\label{prop:small-relu-bisectors}
Suppose finite measures $\mu_1,\ldots,\mu_k$ of positive mass have a
common linear projection $t=u^{\mathsf T}x$ whose image measures have
no atoms. Then a one-hidden-layer ReLU function of width at most $k$
simultaneously bisects the measures.

For a finite sample whose coordinates are pairwise distinct across all
classes, the same width bound gives a strict common bisector whenever
all class sizes are positive and even. Both constructions can be
repeated through depth $H$, with one projection $u$ shared by the whole
tree; in the finite case assume $2^H\mid n_i$ for every $i$.
\end{theorem}

\begin{proof}
For the projected atomless measures, the space
$\operatorname{span}\{1,t,\ldots,t^k\}$ has dimension $k+1$, and every
nonzero polynomial has a zero set of measure zero. The same odd-map
argument supplies a bisecting polynomial $p$ of degree at most $k$.
This is the simultaneous sign-bisection principle associated with
Hobby--Rice~\cite{HobbyRice1965}. For the finite sample choose a
projection separating all coordinates. The Vandermonde determinant
makes its evaluation rows independent, so
Theorem~\ref{thm:finite-polynomial} supplies a strict polynomial bisector.

It remains to realize its signs by a small ReLU function. Let
$r_1<\cdots<r_s$, $s\le k$, be the sign-changing real roots of $p$,
and set $m_j=(r_j+r_{j+1})/2$. There is at least one sign change,
since every class has positive mass and is bisected. Choose $L<r_1$
and $a\in\{-1,1\}$ so that $-a$ is the sign of $p$ to the left of
$r_1$. Define
\begin{equation}\label{eq:relu-common-bisector}
 q(t)=a(L-r_1)+a\relu(t-L)
       +2a\sum_{j=1}^{s-1}(-1)^j\relu(t-m_j).
\end{equation}
The successive nonconstant slopes are $a,-a,a,\ldots$. The first zero
is $r_1$, and reflection of the slope at each midpoint $m_j$ puts the
next zero at $r_{j+1}$. The left constant value also has the required
sign. Thus $q$ has the same signs as $p$ except possibly at roots of
even multiplicity, which have zero measure and are absent from a strict
finite sample. Formula~\eqref{eq:relu-common-bisector} uses $s$ ReLU
units and no additional linear term.

Restriction to a reached region preserves atomlessness of the projected
measures; a subset of the finite sample retains distinct projected
coordinates. This proves the recursive statements.
\end{proof}

In the finite case, ordering the projected points turns their labels
into the colors of a necklace. Its division here is the two-recipient
case of the necklace-splitting theorem~\cite{AlonNecklace}. The bound
of $k$ division points is sharp in this order: if each color occupies
one consecutive block, each block must be cut.

The projection assumption is automatic for absolutely continuous
measures. In fact, a common atomless projection exists for any finite
family of nonatomic finite Borel measures. For almost every direction
$u$, two independently sampled distinct points have different
projections: for each fixed $x\ne y$, the exceptional directions
satisfying $u^{\mathsf T}(x-y)=0$ form a set of spherical measure zero.
Fubini's theorem, together with the zero product measure of the diagonal,
shows that the projected measures have no atoms. In dimension one the
identity projection suffices.

The theorem allows the hidden units to be chosen for the particular
measures or sample. A wide but arbitrarily fixed hidden layer need not
work. For instance, on the two intervals $[-2,-1]$ and $[1,2]$, every
feature $\relu(t+b_j)$ with $b_j>2$ is affine. Every affine combination
of such features is therefore affine on both supports. It cannot bisect
both uniform interval measures, whose medians are different. Choosing
the positions of the changes in slope is essential.

There is also a geometric limitation. The zero set in the construction
is a union of at most $k$ parallel hyperplanes. Its two regions can be
disconnected, and the projection that makes bisection possible may be
poor for local classification. The width bound gives a useful family
of admissible cuts; choosing a useful member of that family remains a
learning problem. The hidden layer may vary from node to node even
when the projection is shared.

\subsection{Gaussian smoothing and finite reconstruction}
\label{sec:certified-construction}

For fixed features, an affine cut can be written
$g(x)=a^{\mathsf T}\phi(x)+b$. One elementary construction searches
over the direction $a$ and then looks for a common threshold.
If class $i$ has $2m_i$ points and its projected values are ordered as
$s_{i,(1)}\le\cdots\le s_{i,(2m_i)}$, the admissible strict thresholds $-b$
are exactly
\[
 \bigcap_{i=1}^k(s_{i,(m_i)},s_{i,(m_i+1)}).
\]
The intersection is nonempty precisely when
$\max_i s_{i,(m_i)}<\min_i s_{i,(m_i+1)}$.
This gives a simple exact test of a proposed direction.

To obtain a continuous problem in the direction as well, replace each
point in feature space by an isotropic Gaussian distribution of scale
$\sigma>0$. Normalize $\|a\|=1$. If $\Phi$ is the standard normal
distribution function, the excess positive mass of class $i$ is
\begin{equation}\label{eq:gaussian-balance-map}
 f_{i,\sigma}(a,b)=
 \frac1{n_i}\sum_{x\in X_i}
       \Phi\left(\frac{a^{\mathsf T}\phi(x)+b}{\sigma}\right)
       -\frac12.
\end{equation}
For an unnormalized nonzero normal, the denominator is
$\sigma\|a\|$. This normalization is part of the mass calculation.
The functions in~\eqref{eq:gaussian-balance-map} are smooth, and their
derivatives are explicit finite sums.

Existence for the smoothed measures follows from the ham-sandwich
theorem in feature space. Equivalently, use homogeneous parameters
$(a,b)\in S^M$: away from the poles the denominator is $\sigma\|a\|$,
and at $a=0$ the excess masses extend continuously to $\pm1/2$.
When $M>k$, one may first restrict the normal to a fixed $k$-dimensional
linear subspace; the corresponding homogeneous parameters form $S^k$.
In this square system a regular homotopy has a one-dimensional zero set.

Our earlier continuation method~\cite{MusinFindingSolutions} joins a
balance map with a known solution to the required map and follows a
curve of zeros. At regular points, tangent prediction and Newton
correction give a numerical approximation to this curve. The relevant
quantitative conditions concern derivatives, the inverse Jacobian and
the path being followed. For fixed feature dimension, evaluating the
map and its first derivatives requires a scan of the sample; the total
work also depends on the number and difficulty of the continuation
steps. The topological existence theorem by itself supplies no bound
for that number.

The smoothed and finite problems are close when few observations lie
near the cut, but small mass discrepancy alone does not give equal
integer counts. The next estimate states the connection directly.

\begin{proposition}[Smooth mass and finite counts]\label{prop:residual-tail}
Let $X_i$ have positive even size $n_i$, and suppose $g(x)\ne0$ at
every point of $X_i$. For $\sigma>0$ put
\[
 r_i=\frac1{n_i}\sum_{x\in X_i}\Phi(g(x)/\sigma)-\frac12,
 \qquad
 t_i=\frac1{n_i}\sum_{x\in X_i}\Phi(-|g(x)|/\sigma).
\]
Then
\begin{equation}\label{eq:residual-count-bound}
 \left|\frac{|\{x\in X_i:g(x)>0\}|}{n_i}-\frac12\right|
       \le |r_i|+t_i.
\end{equation}
In particular, $|r_i|+t_i<1/n_i$ implies exact bisection. If
$\gamma_i=\min_{x\in X_i}|g(x)|$, the simpler sufficient condition is
\begin{equation}\label{eq:residual-margin-certificate}
 |r_i|+\Phi(-\gamma_i/\sigma)<1/n_i.
\end{equation}
\end{proposition}

\begin{proof}
For each nonzero $g(x)$,
\[
 \left|\Phi(g(x)/\sigma)-\mathbf1_{\{g(x)>0\}}\right|
       =\Phi(-|g(x)|/\sigma).
\]
Averaging and applying the triangle inequality proves
\eqref{eq:residual-count-bound}. Since $n_i$ is even, the count
discrepancy is an integer multiple of $1/n_i$, so a smaller absolute
value must be zero. Finally $t_i\le\Phi(-\gamma_i/\sigma)$.
\end{proof}

The tail term measures how much of the smoothed mass can cross the
boundary near observed points. It is indispensable: four positive
scores $\delta,2\delta,3\delta,4\delta$ have a mass discrepancy tending
to zero as $\delta/\sigma\to0$, while all four points remain on the
same side. A direct count of the signs is therefore the natural final
check when exact finite bisection is the goal. The estimate explains
when a successful smoothed calculation already guarantees that check.

A nearby exact cut can also be sought by changing its coefficients.
The essential simplification is to identify points whose signs cannot
change in the chosen neighbourhood. Suppose
\[
 g_\beta(x)=g_0(x)+\sum_{j=1}^r\beta_j f_j(x),\qquad
 |\beta_j|\le h_j,
\]
where the features are fixed, $h_j>0$, and $g_0$ is nonzero on the
whole sample. Let $S$ be the indexed subcollection of observed points
satisfying the displayed inequality:
\[
 E(x)=\sum_jh_j|f_j(x)|,
 \qquad S=\{x\text{ in the sample}:|g_0(x)|\le E(x)\}.
\]
Every point outside $S$ keeps its strict sign throughout the coefficient
box. Let $s_i=|X_i\cap S|$, let $p_i$ be the number of positive points
of class $i$ outside $S$, and put $r_i=n_i/2-p_i$.

\begin{proposition}[Exact reconstruction in a coefficient box]
\label{prop:complete-bounded-repair}
If $r_i\notin\{0,\ldots,s_i\}$ for some $i$, the box contains no
strict common bisector. Otherwise enumerate the
$\prod_i\binom{s_i}{r_i}$ sign assignments on $S$ having $r_i$
positive signs in each class. For an assignment $t:S\to\{-1,1\}$,
solve the linear program
\begin{equation}\label{eq:finite-repair-lp}
 \begin{array}{ll}
 \text{maximize}&\eta\\
 \text{subject to}&t(x)g_\beta(x)\ge\eta\quad(x\in S),\\
 &-h_j\le\beta_j\le h_j\quad(1\le j\le r),\qquad -L\le\eta\le L,
 \end{array}
\end{equation}
where $L=1+\max_{x\in S}(|g_0(x)|+E(x))$ and the empty maximum is
zero. A strict common bisector exists in the box if and only if at
least one of these programs has a positive optimum. If the values
$g_0(x)$, $f_j(x)$ and $h_j$ are rational, the decisions can be verified
by exact rational linear programming.
\end{proposition}

\begin{proof}
Outside $S$, the inequality
$|g_\beta(x)-g_0(x)|\le E(x)<|g_0(x)|$ fixes every sign.
Thus the stated class counts on $S$ are necessary. A positive feasible
$\eta$ gives all signs of an enumerated assignment strictly, and hence
a common bisector. Conversely, a strict bisector in the box has one of
these assignments and a positive minimum signed value on $S$.
The bound on $L$ does not exclude that value. If $S$ is empty, the
same argument says that the unchanged signs either have the required
counts or make them impossible.

Each linear program has a nonempty bounded feasible polytope:
$\beta=0$, $\eta=-L$ is feasible. Its optimum is attained at a vertex.
For rational data, vertices are obtained from rational linear systems;
linear programming duality supplies exact upper bounds as well as
feasible candidates. An optimum at most zero excludes that assignment.
Since the enumeration includes every possible balanced sign assignment,
it also decides nonexistence in the prescribed box.
\end{proof}

This reduction is useful when $S$ is small. It separates a scan of all
observations from the finite geometric problem close to the proposed
cut. The number of assignments can nevertheless be large, and exhausting
a prescribed search budget is different from proving that no cut exists.
Connecting the numerical end of continuation to an effective choice of
such a box is one of the questions discussed in Section~\ref{sec:experiments}.

\paragraph{When the smoothed measure is the intended object.}
\label{sec:statistical-smoothing}
For many applications the observations are samples from a distribution,
or their coordinates are themselves uncertain. Then the smoothed measure
has a statistical meaning and need not be regarded merely as a device
for obtaining exact halves of the observed points. Kernel estimates of
distribution functions and quantiles are classical~\cite{Azzalini1981}.
The distinction between error in the measure and error in its numerical
bisection is also discussed in our earlier work~\cite{MusinFindingSolutions}.

For example, let $K_{ij}$ be a specified conditional distribution for
the unknown location of observation $j$ in class $i$, given the observed
sample. The measure
$Q_i=n_i^{-1}\sum_jK_{ij}$ is the expected empirical distribution of
those locations. A set $A$ with $Q_i(A)=1/2$ contains one half of the
class in conditional expectation. If the locations are conditionally
independent, Hoeffding's inequality gives, simultaneously for all classes
with probability at least $1-\delta$,
\[
 \left|\frac{N_i(A)}{n_i}-\frac12\right|
 \le \left|Q_i(A)-\frac12\right|
       +\sqrt{\frac{\log(2k/\delta)}{2n_i}},
\]
where $N_i(A)$ counts the uncertain locations in $A$~\cite{Hoeffding1963}.
Conditional on the observed sample, $A$ is fixed, so this follows by
applying the inequality to independent indicators and then taking a
union bound. It gives one precise setting in which accurate mass
balance is an appropriate statistical goal.

The distributions $K_{ij}$ must be justified by the measurement model.
In particular, adding Gaussian noise to a recorded coordinate does not
automatically recover the distribution of its unobserved true value.
Similarly, for nonlinear $\phi$, Gaussian smoothing in feature space
defines a different measure from the image of Gaussian uncertainty in
the original coordinates. A statistical tree should restrict the
chosen original class measures at successive nodes, rather than replace
them by new measures at each cut. Exact observed counts and balance of
these fixed measures may lead to different partitions; either objective
can be appropriate once its role is stated.

\subsection{When a local classifier is simpler}\label{sec:local-geometry}

A balanced leaf still contains all classes. The useful question is
whether the decision boundaries needed \emph{inside} it are simpler
than those for the whole sample. Separated graphs give an elementary
setting in which this can be proved.

Write a point of $\R^{r+1}$ as $(u,z)$, where $u\in\R^r$ describes
its position along the graphs and $z$ is the remaining coordinate.
The classes lie near the graphs of functions $f_1<\cdots<f_k$.
Over a small range of $u$, their intermediate surfaces can be replaced
by tangent hyperplanes. The following statement quantifies the required
scale.

\begin{proposition}[Separated local graphs]\label{prop:local-graphs}
Let $f_1<\cdots<f_k$ be $C^2$ functions on a convex domain in $\R^r$,
with $f_{j+1}-f_j\ge\Delta>0$ and
$\|D^2f_j\|_{\mathrm{op}}\le M$.
Suppose class-$i$ points $(u,z)$ satisfy
$\|u-u_0\|\le h$ and $|z-f_i(u)|\le\eta$, with $u$ and $u_0$ in the domain.
If
\begin{equation}\label{eq:local-graph-scale}
 Mh^2+2\eta<\Delta,
\end{equation}
then affine class scores in the original coordinates $(u,z)$ classify
all these points correctly.
\end{proposition}

\begin{proof}
Let $a_j(u)$ be the affine Taylor polynomial at $u_0$ of the midpoint
$(f_j+f_{j+1})/2$. Its error is at most $Mh^2/2$.
For a class-$i$ observation, $z-a_j(u)$ is at least
$\varepsilon=(\Delta-Mh^2-2\eta)/2>0$ when $j<i$, and at most
$-\varepsilon$ when $j\ge i$.
Choose the affine scores
\[
 \ell_1(u,z)=0,\qquad
 \ell_i(u,z)=\sum_{j<i}(z-a_j(u)).
\]
For any $s\ne i$, summing the appropriate increments gives
$\ell_i-\ell_s\ge|i-s|\varepsilon>0$.
\end{proof}

The estimate exhibits the geometric quantities that matter: the gap
between classes, their curvature, the variation allowed within a region,
and the thickness of the observations about the ideal surfaces. With
$\eta=0$, a base radius below $\sqrt{\Delta/M}$ suffices when $M>0$.
Closer or more curved graphs require smaller regions. If the graphs are
affine, curvature imposes no restriction on their base diameter.

In a favourable partition, local networks therefore need only approximate
nearly affine boundaries. They receive the original coordinates and do
not need the unknown functions $f_i$ as extra inputs. The proposition
constructs possible scores from the population geometry; fitting them
from a finite sample is a separate step. Its purpose is to explain why
local approximation can reduce the necessary complexity of a classifier.

Balance does not imply the bound on the base radius in
\eqref{eq:local-graph-scale}. Nor does it prevent one leaf from collecting
several distant pieces on which the useful affine rules differ. This is
particularly relevant for the parallel cuts constructed in
Theorem~\ref{prop:small-relu-bisectors}. An algorithm needs a criterion
for choosing between balanced partitions, and it should account for
the kind of local model that will actually be fitted.

Finite sampling creates another limitation of a geometric explanation.
Suppose one hopes to place each of $R$ separate geometric regions in its
own leaf of a complete balanced tree, with $R=2^H$. Then each such region
must contain exactly $n_i/R$ observed points of class $i$. Even when
the population masses of the regions are equal, independently sampled
counts will usually differ. Some observations must then cross the
proposed geometric grouping. This affects a particular partition, not
the possibility of successful classification: a leaf containing more
than one geometric piece may still admit a simple local classifier.

\subsection{Fitting local models and improving the partition}
\label{sec:joint-learning}

Once a partition is chosen, its local models can be fitted independently.
A familiar difficulty arises near a boundary introduced by the partition:
the local sample approaches the boundary from only one side. Information
from nearby points across it can be useful, as in local regression
\cite{Fan1992}. This motivates fitting with overlap while retaining the
same partition for prediction.

\paragraph{Overlap in the fitting samples.}\label{sec:overlap}
For a leaf region $C_v$, let $X_v$ denote all its training observations,
and enlarge this fitting sample by nearby observations from other leaves.
For instance, for $\delta>0$ one may use
\[
 C_v^{+\delta}=\{x:\dist(x,C_v)<\delta\},\qquad
 J_v=X_v\cup\{x\text{ in the training sample}:x\in C_v^{+\delta}\}.
\]
Balance applies to $X_v$. The enlarged fitting sample $J_v$ may have different class proportions.
A new point is still classified by the single local model for its leaf,
as in~\eqref{eq:cut-network}. Only the data used to estimate that model
have increased; see Figure~\ref{fig:routing-overlap}.

A convenient objective is
\begin{equation}\label{eq:overlap-objective}
 A_v(\theta)+\lambda B_v(\theta),\qquad
 A_v(\theta)=\frac1{|X_v|}\sum_{i=1}^k\sum_{x\in X_{v,i}}
       \ell(h_{v,\theta}(x),i)+\mathcal R(\theta),
\end{equation}
where $B_v$ is the average loss on the added observations and its term
is omitted when there are none. The penalty $\mathcal R$ controls the
size of the fitted coefficients. For class scores $z\in\R^k$, one
standard loss is
$\ell(z,y)=\log\sum_{i=1}^k e^{z_i}-z_y$.
It penalizes giving the observed class a small score relative to the
others. Nothing in the definition of balance singles out this loss;
other fitting criteria can be used.

Overlap trades additional information for a possible change in the local
fit. If the neighbouring data obey a different local rule, putting too
much weight on them can increase error in $C_v$. The following elementary
bound is useful when only the final affine coefficients are being fitted
and the objective is strongly convex.
Here $\mu$-strong convexity of $A$ means that
$A(\theta)-\mu\|\theta\|^2/2$ is convex.

\begin{proposition}[Effect of an added fitting sample]
\label{prop:controlled-overlap}
Let $A:\R^p\to\R$ be differentiable and $\mu$-strongly convex,
$\mu>0$, and let $B:\R^p\to\R$ be differentiable and convex.
Suppose $\theta_0$ minimizes $A$, $\theta_\lambda$ minimizes
$A+\lambda B$, and $\lambda\ge0$. If
$\|\nabla B(\theta_0)\|\le G$, then
\[
 \|\theta_\lambda-\theta_0\|\le\frac{\lambda G}{\mu},\qquad
 0\le A(\theta_\lambda)-A(\theta_0)
       \le\frac{\lambda^2G^2}{2\mu}.
\]
\end{proposition}

\begin{proof}
The gradient of $A+\lambda B$ is $\mu$-strongly monotone. Its values
at $\theta_0$ and $\theta_\lambda$ are
$\lambda\nabla B(\theta_0)$ and zero, respectively. Their scalar
product with $\theta_0-\theta_\lambda$, followed by Cauchy--Schwarz,
gives the first bound. Set $d=\|\theta_\lambda-\theta_0\|$.
Strong convexity of $A+\lambda B$ at its minimizer and the supporting
hyperplane inequality for $B$ give
\[
 A(\theta_\lambda)-A(\theta_0)
       \le\lambda Gd-\frac\mu2d^2
       \le\frac{\lambda^2G^2}{2\mu}.
\]
The lower bound follows from the definition of $\theta_0$.
\end{proof}

For an affine binary score $s_\theta(x)=\theta^{\mathsf T}z(x)$, this
also preserves the prediction at a point whenever
$\lambda G\|z(x)\|/\mu<|s_{\theta_0}(x)|$.
Such a bound is a way to control the influence of neighbouring points
on a fixed convex fitting problem. In a model with a common local rule,
extra observations can instead reduce variance, particularly near a
boundary. Which effect dominates depends on the data and the local
model. Section~\ref{sec:experiments} discusses both behaviours in the
experiments.

\paragraph{Choosing the cuts together with the local models.}
One can construct all cuts first and then fit the local classifiers.
Alternatively, one can alternate between proposing cuts and fitting
local models. In the latter approach, each proposed partition is
checked against~\eqref{eq:strict-balance}, and admissible candidates are
compared by the quality of the local fits. A small affine fit is cheap
to compute; a temporary neural network may better reflect the difficulty
faced by the eventual local classifier. The cost of making this choice
belongs to construction, even though the temporary models are not used
to answer queries.

A common mathematical formulation makes both choices explicit. If $T$
is an admissible balanced tree and $h=(h_v)$ is its family of local
networks, define
\begin{equation}\label{eq:joint-objective}
 \mathcal L(T,h)=\frac1N\sum_{i=1}^k\sum_{x\in X_i}
       \ell(h_{v_T(x)}(x),i),\qquad
 J_\lambda(T,h)=\mathcal L(T,h)+\lambda C(T,h),
 \quad\lambda\ge0.
\end{equation}
Here $C(T,h)$ is a specified measure of prediction cost, including the
cuts and the chosen local architectures. One seeks a small value of
$J_\lambda$ over prescribed families of trees and networks. The tree
with no cuts is allowed and gives an ordinary global classifier. Thus
the formulation permits the data and the cost criterion to determine
whether subdivision is useful. It describes a joint design problem;
the numerical comparisons below use the particular fitting and model
selection procedures stated there.

There is a simple reason why direct differentiation can be ineffective.
At a strictly divided finite sample, sufficiently small changes in the
cut coefficients leave every point on the same side. With the local
models fixed, the classification loss is then constant as a function
of the cut coefficients. Even its infimum over the same local model
families is unchanged in that neighbourhood. Moving between genuinely
different partitions therefore requires crossing sample points, or
using a different objective to propose the changes.

A continuous alternative replaces the choice of one branch by weights
depending smoothly on the cut score. This idea is used in mixtures of
experts and soft decision trees
\cite{JacobsEtAl1991,JordanJacobs1994,IrsoyYildizAlpaydin2012}.
It makes an objective differentiable, but averages several local outputs
during fitting. Since the final classifier in~\eqref{eq:cut-network}
uses one local output, the relation between the two predictions must
be examined. The following bound gives a direct condition.

\begin{proposition}[Replacing an average by one local classifier]
\label{prop:soft-hard-certificate}
At a fixed input, let $z_v\in\R^k$ be the local class scores and
$p_v\ge0$, $\sum_vp_v=1$, their averaging weights. Let $h$ be the
leaf chosen by the signs of the cuts and set
\[
 \bar z=\sum_vp_vz_v,\qquad
 R=\sum_{v\ne h}p_v\operatorname{osc}(z_v-z_h),
 \qquad \operatorname{osc}(a)=\max_i a_i-\min_i a_i.
\]
Then $|\ell(\bar z,y)-\ell(z_h,y)|\le R$ for every label $y$.
If the difference between the largest and second-largest scores of
either $\bar z$ or $z_h$ is greater than $R$, both vectors predict the
same class.
\end{proposition}

\begin{proof}
Put $a=\bar z-z_h$. The difference of the two log-sum-exp terms lies
between $\min_i a_i$ and $\max_i a_i$. Subtracting $a_y$ bounds the
loss difference by $\operatorname{osc}(a)$. Each pairwise difference
of class scores also changes by at most this quantity. Finally,
\[
 \operatorname{osc}(a)
 =\operatorname{osc}\!\left(\sum_{v\ne h}p_v(z_v-z_h)\right)
 \le\sum_{v\ne h}p_v\operatorname{osc}(z_v-z_h)=R.
\]
\end{proof}

Thus it is useful for the averaging weight to concentrate on the chosen
leaf and for the local models not to disagree excessively. A small
loss for the averaged scores alone gives neither property: a large
correct score from an unchosen model can compensate for an incorrect
score from the chosen one. The final one-leaf classifier is therefore
the object whose accuracy and cost are compared in the experiments.

An alternative elementary safeguard is to retain an unchanged model
when refining a region. Copying a parent classifier to both children
preserves every prediction. For empirical loss on the observations reaching the parent, its mean
is the weighted mean of the two child losses, with weights proportional
to child sample sizes. If the child model families contain the copies,
minimizing separately can only decrease this weighted loss. With
overlap or penalties the fitting objective can change upon subdivision,
and numerical optimization need not attain its minimum. It is one reason to compare increasingly strong global
networks as well as increasingly refined partitions.

\subsection{Cost of prediction and construction}\label{sec:costs}

The distinction between the whole fitted model and the portion used by
one query is central to the method. Suppose a complete depth-$H$ tree,
$H\ge1$,
uses one shared $M$-coordinate feature map. Let $C_\phi$ be the number
of multiplications needed to compute it once, and let $S_\phi$ count
its stored parameters. If each leaf classifier uses $C_h$
multiplications and $S_h$ stored coefficients, the dense arithmetic
counts are
\begin{equation}\label{eq:cost}
 C_{\mathrm{query}}=C_\phi+HM+C_h,\qquad
 S_{\mathrm{total}}=S_\phi+(2^H-1)(M+1)+2^H S_h.
\end{equation}
The first expression follows one path; the second counts the entire
tree. A small prediction cost can therefore coexist with substantial
storage.

For affine $k$-class scores,
$C_h=dk$ and $S_h=(d+1)k$.
A one-hidden-layer ReLU network of width $w$ has
\[
 C_h=(d+k)w,\qquad S_h=(d+1)w+(w+1)k.
\]
For binary classification one may use a single score and threshold
instead of two class scores, replacing $k$ by one in these formulas.
A radial quadratic feature costs $d$ squarings. Shared features need
be computed only once; features specific to a node contribute their
cost at every visited node. If no cut is made, there is no need to
compute unused features.

Multiplication counts are transparent comparisons of mathematical
models. Actual running time also depends on branching, memory access,
vectorization and the implementation of the networks. We therefore
report arithmetic cost and measured time separately when both are
available. All models are compared using their full stored coefficients,
including local networks not visited by a particular query.

Construction has a different accounting. At a fixed depth, the regions
partition the training sample. Hence a complete tree satisfies
\[
 \sum_{v\text{ internal}}|X_v|=NH.
\]
If a construction examines the features a bounded number of times at
each node, these scans cost a constant multiple of $MNH$. Solving the
numerical equations, searching over cuts, and fitting the local models
add to this work. In particular, the exhaustive reconstruction in
Proposition~\ref{prop:complete-bounded-repair} can be expensive even
when the final cut is simple.

Overlap also increases fitting work because an observation may contribute
to more than one local sample. It does not increase the number of local
classifiers evaluated for a query. These distinctions make it possible
to state the practical question precisely: whether the effort of
constructing and storing a balanced partition yields a useful reduction
in prediction cost at a given level of accuracy. The examples in the
next section give preliminary answers for particular geometries and
data sets.

\section{Experiments and future work}\label{sec:experiments}

The purpose of the examples is to understand when a balanced partition
makes classification simpler. There are two questions. Does restriction
to a smaller region permit a simpler local decision rule? If so, does
the resulting tree give an accurate prediction with less computation?
The first is a geometric question; the second also depends on how the
partition and the local networks are found from a finite sample.
We describe several examples that separate these issues. They include
favorable configurations, sensitivity to the sampling procedure, and
cases in which an ordinary global network performs better.

\subsection{How the comparisons are made}\label{sec:experimental-rules}

Each observation consists of its coordinates and its class label.
The algorithms receive these observations, but not the equations,
regions, or parameters used to generate a synthetic example. Both the
global network and the local networks use the original coordinates.
When a quadratic feature is used to construct a cut, its evaluation
is included in the cost of the tree.

We distinguish three roles for the observations. The \emph{training
sample} is used to construct the cuts and fit the networks. A separate
\emph{validation sample} is used to choose among fitted models. Finally,
an \emph{evaluation sample}, not used in those choices, measures the
error of the selected classifier. The samples have these different
roles even when their distributions are the same. Repeating a fit on
one training sample studies the effect of initialization; drawing a
new training sample studies an additional source of variation. We
specify which repetition is used in each comparison.

A network of width $w$ has one hidden ReLU layer with $w$ units and
an affine output for each class. Widths $(w,w)$ denote two hidden
layers. An affine local classifier has no hidden layer. In the later
synthetic experiments, each accepted cut divides every class reaching
it into equal halves. A complete tree of depth $q$ consequently has
the prescribed class counts in all $2^q$ leaves. If the numerical
construction stops before that depth, the incomplete tree is recorded
as such. The earlier binary experiments use nearly equal halves
when a class has odd size, as explained below.

The computation for one query consists of the cuts along its path
and the local classifier at its final leaf. It should be distinguished
from the storage of all local networks and from the work required to
construct and fit the tree. Our tables retain the arithmetic convention
used in the corresponding experiment: multiplication--accumulation
terms for the earlier affine comparisons, and products, including
quadratic features, for the later algebraic comparisons. These counts
are comparable within each table. Measured execution times are given
only for the pixel-color example; they depend on the implementation
and computer used.

\subsection{Plane curves and the effect of overlap}\label{sec:binary-planar}

A simple example consists of two translated graphs,
\[
 C_0=\{(t,\sin mt):0\le t\le2\pi\},\qquad
 C_1=C_0+(0,c),\qquad c>0,
\]
with the same uniform distribution of the parameter $t$.
Successive midpoint cuts in $t$ bisect both class distributions.
These curves illustrate the hypothesis of Section~\ref{sec:introduction}.
Every point of either curve has a point of the other at vertical distance
$c$. As $m$ increases at fixed $c$, the separating graph
$y=\sin(mt)+c/2$ changes direction more often. On sufficiently short
intervals its tangent line separates the two classes, so the midpoint
cuts are compatible with both balance and local affine separation.

For $m=4$ and $c=0.35$, three experimental runs used $8193$
training observations, $4096$ validation observations, and $16384$
evaluation observations per class. These earlier binary experiments
allow a difference of one observation when a reached class has odd
size. The trees selected for accuracy
made no evaluation errors. There is also an explicitly constructed
ordinary ReLU network of width $34$ that separates the entire pair
of curves. This comparison illustrates the use of local computation:
the tree evaluates only the local rule it needs, whereas the ordinary
network evaluates its full representation. It does not establish that
the separating function requires a tree.

The distinction is equally visible in nested planar annuli. For
supports $1\le\|x\|\le2$ and $4\le\|x\|\le5$, the function
\[
 \max(|x_1|,|x_2|)-2.5
\]
separates the classes and has both a small tree representation and
a width-four ReLU representation. With displaced annuli and more
complicated surfaces, small global networks remained competitive.
The benefit of a partition therefore depends on the example and on
the available local and global models.

An artificial cut can introduce a statistical difficulty of its own.
Near the cut, a local fit has observations mainly on one side of the
query. Adding nearby observations from the neighboring region may
reduce this effect. This is the motivation for the overlap construction
of Section~\ref{sec:models}: the regions used for prediction remain
disjoint, while the samples used to fit adjacent local classifiers
may overlap. Boundary effects of this kind are familiar in local
regression~\cite{Fan1992}.

The following calculation makes both the benefit and its limitation
visible. Divide the plane by $x=0$, and take observations
\[
 z=f(x)+0.09c+\epsilon,\qquad
 c\in\{-1,1\},\qquad \epsilon\sim N(0,0.03^2).
\]
In each half-plane there are $n$ observations of each class, at base
coordinates $\pm j/n$, $j=1,\ldots,n$. Each local fit estimates an
affine midpoint and the class offset by least squares. In the overlap
variant it also uses observations from the opposite-side interval
$0<|x|\le1/2$, giving each added observation the same weight as an
original observation. The noise in fitting and prediction is averaged
analytically. An independent query has $x$ uniform on $[-1,1]$;
errors are then averaged conditionally on each reported region.
Table~\ref{tab:overlap-comparison} gives expected errors;
these are a model calculation, not error counts from trained neural
networks.

\begin{table}[htbp]
\centering\small
\begin{tabular}{@{}lr rrrr@{}}
\toprule
Midpoint & $n$ & \multicolumn{2}{c}{Near the cut}
                         & \multicolumn{2}{c}{Away from the cut}\\
 & & Separate & Overlap & Separate & Overlap\\
\midrule
$f(x)=0$ & 4 & 1.07 & 0.22 & 0.33 & 0.24\\
$f(x)=0$ & 16 & 0.24 & 0.15 & 0.17 & 0.16\\
$f(x)=0.3|x|$ & 4 & 1.07 & 43.51 & 0.33 & 5.55\\
$f(x)=0.3|x|$ & 16 & 0.24 & 28.39 & 0.17 & 2.89\\
\bottomrule
\end{tabular}
\caption{Expected classification error in percent for separate and
overlapping local fits. The near-cut interval is $|x|<0.05$; the
away-from-cut region is $0.2\le|x|\le1$. Values are rounded to two
decimal places.}
\label{tab:overlap-comparison}
\end{table}

For a common affine midpoint, the added observations reduce the
uncertainty of the fit, especially near the boundary. For
$f(x)=0.3|x|$, each side has a correct affine model, but these two
models have different slopes. Extending a fit too far across the
boundary introduces a large bias. Overlap must therefore reflect
compatibility of the neighboring geometry. Its size and weight are
part of the learning problem. A separate experiment with fitted
ReLU features gave only a two-error improvement among almost
error-free comparisons, with all errors away from the cut. Thus the
clear effect in Table~\ref{tab:overlap-comparison} remains to be
established more broadly for trained neural classifiers.

\subsection{Several classes with different local rules}
\label{sec:regime-experiments}

To make local simplicity compatible with global complexity, place
$16$ or $64$ separated regions on a grid in $\R^3$. In region $r$,
the observations of class $j$ have the form
\[
 x=c_r+Q_r(u,v,a_{\pi_r(j)})^{\mathsf T},\qquad
 a=(-0.12,0,0.12),
\]
where $Q_r$ is a rotation, $\pi_r$ is a permutation of the classes,
and $u,v$ are uniform on $[-0.35,0.35]$. Every region contains all
three classes. Within one region the classification is affine, but
the normal direction and class order vary from region to region.

The initial construction used $4096$ training observations per class,
equally distributed among the regions. Validation and evaluation
observations were sampled with independent uniform choices of their
regions. These equal training counts allow a balanced tree to separate
whole regions. The fitting procedure receives only coordinates and
labels; it is not given the regions or rotations. Three repetitions
of the fitting procedure give Table~\ref{tab:local-regions}.

\begin{table}[htbp]
\centering\small
\begin{tabular}{@{}llrr@{}}
\toprule
Regions & Model & Evaluation errors & Terms per query\\
\midrule
$16$ & Selected global network & $7/0/0$ & $768$--$4480$\\
 & Depth-4 tree, affine local classifiers & $0/0/0$ & $21$\\
$64$ & Selected global network & $8668/9245/7705$ & $17152$\\
 & Depth-6 tree, affine local classifiers & $0/0/0$ & $27$\\
\bottomrule
\end{tabular}
\caption{Designed three-class examples. Each population has one
evaluation sample of $98304$ observations, shared by all fitted models.
Slash-separated entries refer to three fitting repetitions. The
selected global architecture varies in the 16-region example, hence
the range of arithmetic costs. For comparability we count dense affine
evaluation throughout: three products per cut and nine for the three
leaf scores, giving $3H+9$ products at depth $H$.}
\label{tab:local-regions}
\end{table}

Here the partition finds precisely the geometric simplification for
which it was intended. For the fitted affine trees, checking linear
inequalities at the corners of every generating patch also verifies
correct classification throughout those patches. The result extends
beyond the finitely many evaluation observations. This is an
especially favorable construction, however, because equal counts
in every class--region pair make complete region separation possible.

A subsequent experiment used a new 16-region environment and $512$
training observations per class. It compared equal region counts
with independent region choices, while retaining exact class balance
at every accepted cut. With independent choices, the selected trees
made $8662$, $11973$, and $6962$ errors on $49152$ new observations
in the three repetitions: error rates between $14.16\%$ and $24.36\%$.
The corresponding global counts were $12595$, $10349$, and $14077$.
The trees used substantially fewer products per query, but their
errors were too large to make this a satisfactory classification
example. Equal total class counts do not force equal counts within
each geometric region; requiring exact halves may split a region
and combine pieces with different local rules.

The cut-selection procedure adds another difficulty. In a later
comparison on independently sampled regions and shifted ellipsoidal
layers, only six of twelve attempts reached the requested depth.
Choosing cuts by the fit of temporary nonlinear local networks did
not consistently improve either completion or validation error.
On the ellipsoidal population an ordinary global ReLU network of
width $18$ can be constructed explicitly to separate all three
supports. This representation was not supplied to the fitting
algorithms. It shows why a disappointing fitted global network
should not be interpreted as a lower bound on its expressive power.

\subsection{Three planar classes and algebraic cuts}
\label{sec:planar-algebraic}

\begin{figure}[htbp]
\centering
\includegraphics[width=0.94\linewidth]{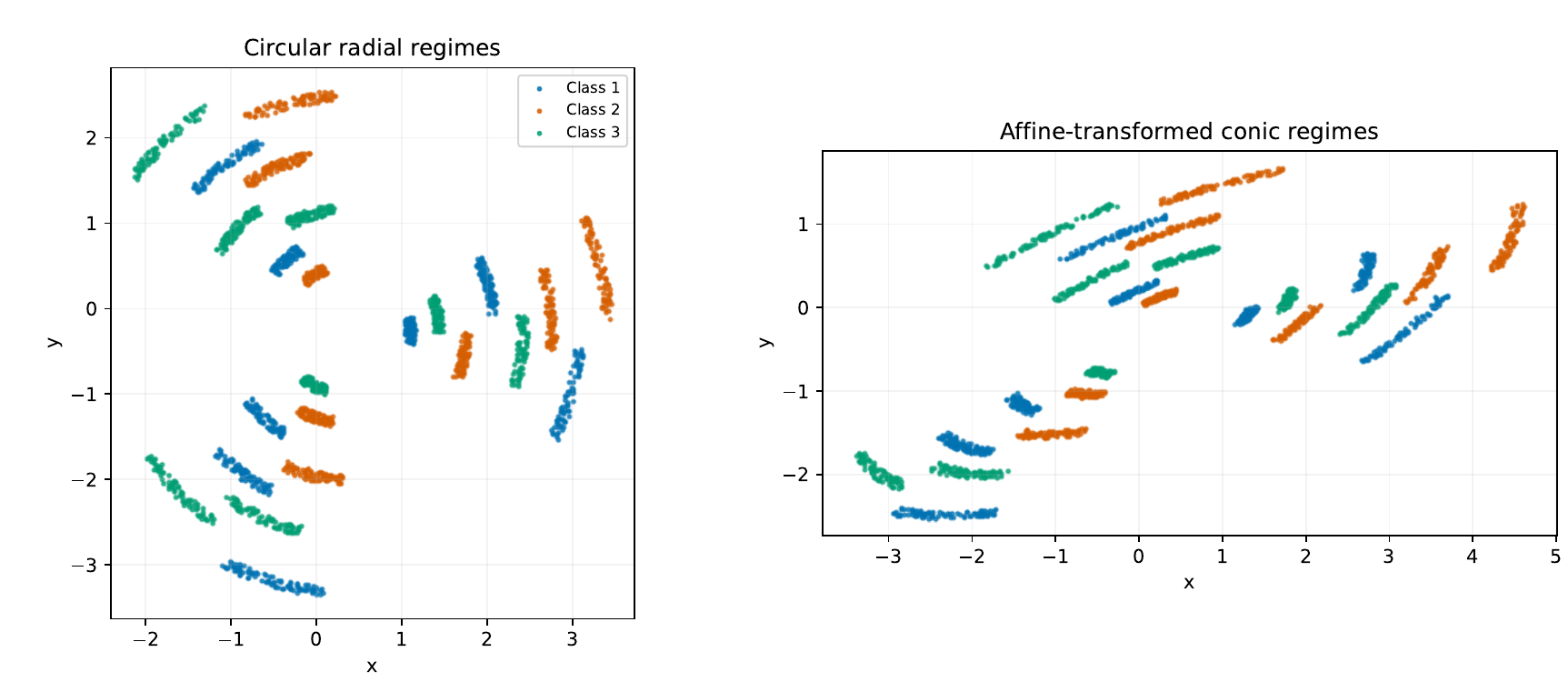}
\caption{Three planar classes in eight radial regions, before and
after an affine change of coordinates. Each region contains all
three classes. Circular or quadratic cuts can separate the regions
while bisecting every class.}
\label{fig:planar-distributions}
\end{figure}

The next example uses positive-area regions in the plane. Eight
thin annular regions have radii $r_i=0.75+0.33i$, $0\le i\le7$.
Within each annular region, each class occupies one of three
separated angular sectors. The angular position and the class order
change with the radius. Each class assigns equal probability to the
eight regions, and sampling within a sector is uniform in area.
An invertible affine change of coordinates gives the second
population in Figure~\ref{fig:planar-distributions}.

In the first experiment there were $2048$ training observations per
class, with $256$ in each radial region. A depth-3 tree then has
eight leaves with the same prescribed counts. Circular cuts gave
zero errors in both fitting repetitions on $98304$ evaluation
observations, using affine local classifiers and $17$ products per
query. For the affine-transformed population, a learned quadratic
form shared by all cuts gave zero errors with $21$ products.
A selected global network with two hidden layers of width $32$
also gave zero or two errors, at $1184$ products. The substantial
difference in this example is the cost of evaluating the classifier.

Even here, finding a balanced cut does not determine which cut to
use. Another quadratic-tree construction gave $0$ errors in one
fit and $9126$ in the other. The latter tree had a slightly larger
separation between its training points and the first cut, yet its
leaves combined incompatible class configurations. In one such leaf,
crossing-segment inequalities show that no affine classifier can
separate all its training observations. Replacing the affine local
classifiers on these same cuts by width-four ReLU networks gave
zero evaluation errors. Thus the choice of a cut and the capacity
of the intended local model must be considered together.

We also studied six independently drawn training samples from the
affine-transformed population, without equal counts in the radial
regions. Each sample again had $2048$ observations per class. Twelve
candidate quadratic forms were obtained from the observed data;
each completed depth-3 tree used width-four local networks. Two
selection rules were specified before evaluation. The first chose
a form using an internal division of the training sample, then
rebuilt the tree on the whole sample. The second chose a completed
tree by its fitted loss on the whole training sample. An unsuccessful
rebuilding was retained, without substitution by a different form.

\begin{table}[htbp]
\centering\small
\begin{tabular}{@{}rrrrr@{}}
\toprule
Sample & \shortstack{First selection\\rule}
 & \shortstack{Second selection\\rule}
 & \shortstack{Global network\\width $8$}
 & \shortstack{Global network\\width $32$}\\
\midrule
1 & $562$ & $3$ & $22962$ & $59$\\
2 & $3249$ & $430$ & $24928$ & $12$\\
3 & $4340^{\dagger}$ & $292$ & $16188$ & $29$\\
4 & $13696^{\dagger}$ & $46$ & $20114$ & $21$\\
5 & $37371^{\dagger}$ & $147$ & $20727$ & $203$\\
6 & $38063^{\dagger}$ & $205$ & $15729$ & $39$\\
\midrule
Products & varies & $35$ & $40$ & $160$\\
Stored coefficients & varies & $247$ & $51$ & $195$\\
\bottomrule
\end{tabular}
\caption{Errors for six independent training samples, measured on
one common evaluation sample of $98304$ observations. A dagger
denotes an incomplete tree after rebuilding. All trees selected by
the second rule reached depth three.}
\label{tab:independent-planar}
\end{table}

The completed trees chosen by the second rule make fewer errors
than the width-eight global networks at a comparable number of
products. They store more coefficients. The width-32 global networks
are usually more accurate than the trees and also require less
storage, although their evaluation uses more products.
The first selection rule exposes a different problem: a cut
construction that succeeds on a smaller sample may stop when it
is repeated on a larger one. These outcomes motivate improvements
to the finite construction as well as to the statistical choice of
the quadratic form.

\subsection{Examples from observed data}\label{sec:observed-data}

The Skin Segmentation data contain $245057$ labeled pixels with
three color coordinates~\cite{BhattDhall2009}. To prevent identical
inputs from occurring in both fitting and evaluation samples, all
occurrences of an exact color triple were placed together. The
resulting training, validation, and evaluation samples contained
$149176$, $45417$, and $50464$ pixels. Original multiplicities and
labels were retained.

The selected global networks had two hidden layers of width $32$.
The selected trees had $9$, $10$, and $11$ leaves in the three
fitting repetitions, with maximum depth four and local widths
$4$, $4$, and $8$. These are earlier, nearly balanced trees.
Every even class count reaching a cut is divided exactly in half;
the only unequal division is $3859=1929+1930$. Their leaves occur
at different depths, so they are not the complete balanced trees
used in the preceding synthetic examples.

\begin{table}[htbp]
\centering\small
\begin{tabular}{@{}rrrrr@{}}
\toprule
Fit & Global errors & Tree errors
 & \shortstack{Global/tree\\terms per query}
 & \shortstack{Global/tree\\time (ns)}\\
\midrule
1 & $45$ & $29$ & $1152/28$ & $276.3/29.5$\\
2 & $37$ & $30$ & $1152/28$ & $277.8/34.0$\\
3 & $43$ & $17$ & $1152/44$ & $279.6/39.4$\\
\bottomrule
\end{tabular}
\caption{Pixel-color classification on the same $50464$ evaluation
pixels. Arithmetic is the maximum-path count. Times are medians
for evaluating one observation in a common compiled implementation,
including normalization and the cuts.}
\label{tab:skin}
\end{table}

In these comparisons the trees have both smaller error counts and
shorter evaluation times than the selected wide global networks.
Smaller global networks are faster still, but have larger errors:
width-four networks make $78/87/114$ errors and width-eight
networks make $73/69/77$. Thus the useful comparison concerns
the relation between accuracy and computation. The timing results
refer to one Linux computer and do not include application overhead.
The data do not identify the original photographs or people, so
the experiment cannot measure transfer to new photographs or
subjects. Pixels with repeated colors are also not independent
statistical observations.

A less favorable example is the ten-class Pendigits collection,
whose observations have $16$ coordinates~\cite{AlimogluAlpaydin1996}.
A single quadratic cut with two width-16 local networks was
selected using development observations. It had $36$ errors among
$1499$ validation observations, compared with $48$ for the
selected global network of similar arithmetic cost. Both used
the same $5968$ training observations. The subsequent evaluation
on the official independent sample reversed this comparison.

\begin{table}[htbp]
\centering\small
\begin{tabular}{@{}lrrr@{}}
\toprule
Classifier & Errors & Products & Stored coefficients\\
\midrule
One quadratic cut, width-16 local networks & $133$ & $465$ & $918$\\
Global network of width $17$ & $105$ & $442$ & $469$\\
Global network of width $32$ & $92$ & $832$ & $874$\\
\bottomrule
\end{tabular}
\caption{Pendigits classification on $3498$ evaluation observations.
The quadratic tree's validation advantage does not persist on this sample.}
\label{tab:pendigits}
\end{table}

The gas-sensor examples are also incomplete as evidence for the
method. In a six-class collection with measurements taken over
time~\cite{Fonollosa2015GasData,RodriguezLujan2014Calibration,Vergara2012Drift},
models were fitted on earlier acquisition groups and compared on
later groups. Across three fitting repetitions, selected global
networks had validation errors of $10.67$--$11.98\%$, while the
best algebraic-tree family in each repetition had
$13.41$--$16.76\%$. The final two acquisition groups have not been
used for evaluation. These results direct attention to changing
distributions and to stronger local models before a final comparison
is attempted.

\subsection{Questions for further work}\label{sec:future-work}

The examples suggest a research program with several related parts.
Evaluating a single leaf saves work in any tree of local models.
To isolate the effect of class balance, one could keep the cut family,
depth, local networks, and fitting budget fixed and compare four rules:
exact class bisection, class bisection within a prescribed tolerance,
bisection of the total count, and selection by local fitting loss alone.
The reported comparisons with global networks do not isolate this effect.
The central issue is to turn the existence of balanced partitions
into a reliable method for finding partitions with useful local
geometry. The following questions concern that passage.

\paragraph{When does proximity favor a balanced partition?}
Can one give geometric conditions under which a balanced tree with
affine leaf classifiers achieves a prescribed error at lower prediction
cost than a global ReLU network of specified depth? A theory would
have to relate the fraction of each class near the other to the
variation of the local separating directions and to the regions that
balanced cuts can produce. For an experimental test, the translated
curves in Section~\ref{sec:binary-planar} allow the gap $c$ and the
frequency $m$ to vary independently. Broader comparisons should also
include surfaces and independently drawn samples, and measure both
prediction cost and error on new observations.

\paragraph{Choosing useful bisectors.}
A cut may have excellent balance and a large separation from the
training points while leaving a difficult problem in one child.
A useful selection rule should take account of the local classifier
to be fitted and of the subdivisions still required below that
child. Limited look-ahead and the possibility of returning to an
earlier choice are natural alternatives to accepting the first
apparently favorable cut. The incomplete trees above provide
concrete examples on which to study these questions. For smooth
populations in general position, theoretical existence may persist
through every subdivision even when the bounded numerical search
stops.

\paragraph{Connecting smoothing with exact finite division.}
Gaussian smoothing gives a continuous equation that can be followed
numerically. An additional step is required when the intended
partition must divide the original finite classes exactly. A study
of $24$ configurations reached a numerical endpoint in $19$ cases
and returned a strictly balanced finite cut in $12$. The examples
included three-class spatial clouds and planar lifts for three or
five classes. Successful finite corrections increased the discrepancy
of the fixed smoothed masses in ten of those twelve cases. This
illustrates the distinction between the two objectives. A practical
method should choose the desired notion of balance, control the
correction needed for finite counts, and record the cost of exact
verification. Exhaustive search in a prescribed coefficient box
is useful when few points have uncertain signs; its general cost
can be large.

\paragraph{Learning a partition and its local classifiers together.}
Preliminary cuts are convenient, but the local classification error
may suggest moving them. This leads to optimization under the
balance constraints. Replacing the selected local output by a
weighted average of several local outputs produces a different
predictor. A controlled comparison on translated sine curves gave
zero error with direct local fitting, mean error $0.174\%$ after
fitting a mixture and then selecting only one local output, and
zero error when the fitting objective accounted for that final
selection. This motivates objectives that correspond directly to
the deployed classifier. Coordinate scaling, penalties, and
initialization also require attention: in a separate fixed-width
experiment, changing the optimization step induced by normalization
reduced validation error counts from $(100,189,234)$ to $(0,9,1)$,
among $3072$ observations in each comparison.
Such sensitivity can be as important as the architectural choice.

\paragraph{Selecting overlap and measuring the cost.}
Neighboring observations should be borrowed where their local rules
are compatible, with weights that reflect their number and relevance.
The positive and negative cases in Table~\ref{tab:overlap-comparison}
make this a geometric and statistical problem rather than a general
prescription to enlarge every local sample. It is also necessary to
measure the complete cost: constructing the tree, fitting all local
networks, storing them, and answering queries. A saving in products
along one path can coexist with larger storage or longer preparation.
Comparisons should include well-trained global networks, explicit
simple formulas when available, and measured execution times under
the same conditions.

\paragraph{Independent evaluation and reproducibility.}
The next empirical step is evaluation on new geometric configurations
and independently acquired groups of observations. Pixel-level or
random-row divisions alone cannot answer questions about transfer
to new people, images, sensors, or periods of time. The present
examples have preserved their successful and unsuccessful fits and
the numerical procedures used to obtain them. Preparing selected
code and data for a public release would permit independent
reproduction and comparison. Until such a release is arranged,
the examples here should be read as a mathematical exploration
supported by preliminary computations. They identify situations
where balanced local classification is useful and specific obstacles
that a more general method must overcome.

\FloatBarrier
\paragraph{Use of AI.} The author used tools from OpenAI for assistance with the computations, the exposition of the proofs, and the preparation of the manuscript. The author is responsible for the arguments and the final text.

\begingroup\small\setlength{\bibsep}{5pt}

\endgroup

\bigskip
\noindent\textbf{Oleg R. Musin}\\
The University of Texas Rio Grande Valley, Brownsville, Texas, USA\\
\href{mailto:oleg.musin@utrgv.edu}{\texttt{oleg.musin@utrgv.edu}}
\end{document}